\documentclass[10pt]{amsart}
\usepackage{amssymb,amsmath,amsthm}
\usepackage{mathrsfs,dsfont,a4wide}

\catcode`@=11 
\def\blfootnote{\xdef\@thefnmark{}\@footnotetext}
\catcode`@=12 

\theoremstyle{plain}
\begingroup
\newtheorem{theorem}{Theorem}[section]
\newtheorem{lemma}[theorem]{Lemma}
\newtheorem{proposition}[theorem]{Proposition}

\endgroup

\theoremstyle{definition}
\begingroup
\newtheorem{definition}[theorem]{Definition}

\endgroup

\numberwithin{equation}{section}

\newenvironment{dimo}{\noindent {\it Proof.}}%
{\hfill \hfill \qued\par\medskip\indent}
\newcommand{\qued}{~\vrule height6pt depth0pt width6pt}

\newcommand{\as}{{\mathcal A}}
\newcommand{\hs}{{\mathcal H}}
\newcommand{\ks}{{\mathcal K}}

\newcommand{\gs}{{\mathcal G}}
\newcommand{\fs}{{\mathcal F}}
\newcommand{\leb}{{\mathcal L}}

\newcommand{\Es}{{\mathcal E}}
\newcommand{\vub}{{\mathcal V}}

\newcommand{\rs}{{\mathcal R}}
\newcommand{\ws}{{\mathcal W}}

\newcommand{\R}{{\mathbb R}}
\newcommand{\N}{{\mathbb N}}

\newcommand{\msd}{\R^{m{\times}n}}

\newcommand{\Om}{\Omega}

\newcommand{\OmBb}{\overline{\Om}_B}

\newcommand{\weak}{\rightharpoonup}

\mathchardef\emptyset="001F

\newcommand{\tsub}{\mathrel{\widetilde{\subset}}}
\newcommand{\teq}{\mathrel{\widetilde{=}}}

\newcommand{\hn}{\hs^{n-1}}

\newcommand{\Eb}{{\mathcal E}^{el}}
\newcommand{\Esup}{{\mathcal E}^s}

\newcommand{\psc}[2]{\langle #1,#2 \rangle}

\newcommand{\Sg}[2]{S^{#1}(#2)}

\newcommand{\radm}{\rs(\OmBb)}

\title
[A variational model for quasistatic crack growth]
{A variational model for quasistatic crack growth
\\
in nonlinear elasticity:
\\
some qualitative properties of the solutions}
\author[G. Dal Maso]
{Gianni Dal Maso}
\address[Gianni Dal Maso]{S.I.S.S.A., Via Beirut 2-4, 34014, Trieste,
Italy}
\email[G. Dal Maso]{dalmaso@sissa.it}
\author[A. Giacomini]
{Alessandro Giacomini}
\address[Alessandro Giacomini]{Dipartimento di Matematica, Universit\`a degli Studi di  Brescia, via Valotti 9, 25133 Brescia,
Italy}
\email[A. Giacomini]{alessandro.giacomini@ing.unibs.it}
\author[M. Ponsiglione]
{Marcello Ponsiglione}
\address[Marcello Ponsiglione]{Dipartimento di Matematica, Universit\`a di Roma ``La Sapienza'', Piazzale A.~Moro 2, 00185 Roma,
Italy}
\email[M. Ponsiglione]{ponsigli@mat.uniroma1.it}
\begin{document}
\vskip .2truecm
\begin{abstract}
\small{
We present the main existence result for quasistatic crack growth in the model proposed by Dal Maso, Francfort, and Toader, and prove some qualitative properties of the solutions.
}
\end{abstract}
\maketitle
{\small \tableofcontents}

\section{Introduction}

\blfootnote{Preprint SISSA 41/2008/M (June 2008)}
In this paper we present an existence result, proved in  \cite{DMFT}, for an evolution problem 
in fracture mechanics, and study some very weak regularity properties of the solutions. 
The mathematical formulation of the problem is based on a variational model for quasistatic 
crack growth developed by Francfort and Marigo \cite{FM}.
The model is based on Griffith's idea \cite{Gri} that the equilibrium 
of a crack is determined by the
competition between the elastic energy released
if the crack grows and the energy dissipated to produce a new portion of crack.
This model not only predicts the crack growth along its path, but also  determines 
the crack path on the basis of an energy criterion, 
and can also be used to study the 
process of crack initiation (see \cite{ChGP}).

The first mathematical results on this model were obtained in \cite{DMT} for linear elasticity 
in the antiplane case in dimension two, assuming an a priori bound on the number of 
connected components of the crack set. This simplifies the mathematical treatment of the 
problem, but has no mechanical justification. 
These results were extended by Chambolle 
\cite{Ch} to the case of plane elasticity. A remarkable improvement was obtained by 
Francfort and Larsen \cite{FL}, who developed a weak formulation in the space 
$SBV(\Om)$ of special functions with bounded variation, introduced by De Giorgi 
and Ambrosio \cite{DGA} to study a wide class of free discontinuity problems. 
This new formulation allows to study the problem in any space dimension, and 
without any restriction on the number of connected components of the crack sets. 
The results of  \cite{FL} deal only with a scalar displacement field, since the
 compactness theorem in $SBV(\Om)$ requires a bound in $L^\infty(\Om)$ 
 that is obtained through a truncation argument. Moreover, the techniques 
 used in all papers considered so far cannot be extended to the case where 
 exterior volume or surface forces act on the body. 

The result we present here allows to deal with the case of a vector valued 
deformation field defined in $\Om\subset\R^n$ with values in $\R^m$,
 including both the case of antiplane shear ($n=2$ and $m=1$) and $n$ 
 dimensional elasticity ($n=m\ge 2$). The bulk energy is not necessarily 
 quadratic, although the polynomial growth \eqref{scisandb2} excludes, 
 for the moment, the case of genuine finite elasticity. Moreover body and 
 surface forces can be considered, under some natural assumptions 
 which imply that the deformed body remains in a bounded region 
 even if the cracks split it into several connected components. 

The variational formulation of the problem presented in this paper
is based on the space $GSBV(\Om;\R^m)$ introduced by De Giorgi and Ambrosio \cite{DGA}, and fits the general framework of rate 
independent evolution problems developed by Mielke~\cite{M}.

The second part of the paper is devoted to the proof of some
 qualitative properties of the solutions, that are not studied in \cite{DMFT}.
 We show that, given any 
 quasistatic evolution, we can construct another quasistatic evolution, 
 whose crack sets are left continuous in time. We also prove that 
 the crack set is minimal, in the sense that, at each time $t$, it is 
 the smallest set containing all discontinuity sets of the 
 deformation at times $s\le t$. Moreover, we prove that, 
 given any quasistatic evolution, we can construct another 
 quasistatic evolution, with the same crack set, whose 
 deformation is measurable in time. Finally, under some 
 convexity assumptions, we prove that, if the crack set is left 
 continuous, so are the deformation, the deformation gradient, and the stress.

\section{The quasistatic crack growth}\label{qsedmft}

In this section we describe in detail the hypotheses of the quasistatic evolution problem studied in  \cite{DMFT} and state, without proof,  the main existence result of~\cite{DMFT}.

\par\medskip\noindent
{\bf The reference configuration.} The reference configuration is a bounded open subset 
$\Om$ of $\R^n$ with Lipschitz boundary. Let $\partial_N \Om \subset \partial \Om$
be closed
and let $\partial_D\Om:= \partial\Om\setminus \partial_N\Om$. We fix 
an open subset $\Om_B$ of $\Om$ with Lipschitz boundary and a closed set
$\partial_S\Om \subset \partial_N\Om$  such that
$\OmBb\cap \partial_S\Om = \emptyset$.
The set
$\overline\Om_B$ represents the brittle part of the body, $\partial_D \Om$ is
the part of the boundary where
the deformation is prescribed, while the surface forces are applied to $\partial_S \Om$.

\par\medskip\noindent
{\bf Admissible cracks.}
The set of admissible cracks is given by
\begin{equation*}
\radm:= \{\Gamma:\Gamma\text{ is rectifiable},\,
\Gamma \tsub \OmBb,
\, \hn(\Gamma)<+\infty\}.
\end{equation*}
Here and henceforth $\hn$ is the $(n-1)$-dimensional Hausdorff measure, 
$\tsub$ means inclusion up to a set of
$\hn$-measure zero, while rectifiable
means that there exists a sequence $(M_i)$ of $C^1$ manifolds of dimension $(n-1)$
such that $\Gamma \tsub \bigcup_i M_i$.
If $\Gamma$ is rectifiable, we can define $\hn$-almost everywhere on $\Gamma$ a Borel measurable
unit normal vector field 
$\nu$ (see, e.g., \cite[Definition~2.86]{AFP}), which is unique up to a pointwise choice of the orientation.

\par\medskip\noindent
{\bf Admissible deformations.}  Given a crack $\Gamma$, an admissible deformation is
given by any function 
$u \in GSBV(\Om;\R^m)$ such that $S(u) \tsub \Gamma$.
We refer to \cite[Chapter~4]{AFP} for the definitions and properties of the spaces 
$SBV(\Om;\R^m)$ and
$GSBV(\Om;\R^m)$, as well as for the definition of the jump set $S(u)$, of the approximate gradient $\nabla u$, and of the trace on $\partial\Om$ of a function $u\in GSBV(\Om;\R^m)$.
Let us fix $p>1$ and $q>1$. We define
\begin{eqnarray*}
&
GSBV^p(\Om;\R^m):=
\{u\in GSBV(\Om;\R^m)\,:\, \nabla u \in L^p(\Om;\R^{m{\times} n}),\,
\hs^{n-1}(S(u)) <+\infty\},
\\
& GSBV^p_q(\Om;\R^m):=GSBV^p(\Om;\R^m) \cap L^q(\Om;\R^m).
\end{eqnarray*}
We say that $u_k \weak u$ weakly in $GSBV^p_q(\Om;\R^m)$ if
\begin{align}
\nonumber
u_k \to u &\quad\text{in measure on } \Om,\\
\label{1.2}
\nabla u_k \weak \nabla u &\quad\text{weakly in }
L^p(\Om; \R^{m {\times} n}), \\
\nonumber
u_k \weak u &\quad\text{weakly in }L^q(\Om;\R^m).
\end{align}

\par\medskip\noindent
{\bf The surface energy.}
The energy spent to produce a crack $\Gamma$ is given by
\begin{equation}
\label{crackener}
\Esup(\Gamma):=\int_{\Gamma\setminus\partial_N\Om} \kappa(x,\nu(x)) \,d\hn(x),
\end{equation}
where $\nu$ is a unit normal vector field on $\Gamma$.
Here $\kappa(x,\nu(x))$ represents the {\it toughness\/} of the material, which depends
on the position $x$ and on the tangent space to the crack, determined by $\nu(x)$.
Since we are dealing only with {\it brittle\/} cracks, the toughness does not depend
on the size of the jump of~$u$.

\par
We assume that
$\kappa\colon\OmBb{\times}\R^n \to \R$ 
is continuous, that
$\kappa(x,\cdot)$ is a norm in $\R^n$
for all $x \in \OmBb$, and that
\begin{equation}
\label{ik3}
K_1|\nu| \le \kappa(x, \nu) \le K_2|\nu| \quad\hbox{for all
}x\in  \OmBb\hbox{
and }\nu\in \R^n,
\end{equation}
with $K_1$, $K_2>0$. Notice that, since $\kappa$ is even in the
second variable, the integral \eqref{crackener} depends only on the
geometry of $\Gamma$, and is independent of the choice of the orientation of $\nu(x)$.

\par\medskip\noindent
{\bf The bulk energy.}
Let $p>1$ be fixed.
Given a deformation $u \in GSBV^p(\Om;\R^m)$ the
associated {\it bulk energy}
is given by
\begin{equation}
\label{bws}
\ws(\nabla u):= \int_\Om W(x,\nabla u(x)) \, dx,
\end{equation}
where $W\colon \Om {\times} \msd \to [0,+\infty)$ is a
Carath\'eodory function satisfying the following conditions:
\begin{eqnarray}
\label{scisandb1}
& 
W(x,\cdot) \text{ is quasiconvex and }
C^1 \text{ on }\msd\text{ for every } x \in \Om,
\vspace{2pt}
\\
\label{scisandb2}
&
a_0^W|\xi|^p - b_0^W(x)\le W(x,\xi)\le
a_1^W |\xi|^p + b_1^W(x) \text{ for every } (x,\xi) \in \Om{\times} \msd.
\end{eqnarray}
Here $a_0^W>0$ and $a_1^W>0$ are constants, while
$b_0^W$ and $b_1^W$
are nonnegative functions in $L^1(\Om)$. The quasiconvexity assumption means that
$$
W(x,\xi) \le \int_\Om W(x,\xi +\nabla \varphi(y))\,dy
$$
 for all $x\in\Om$, $\xi \in \msd$, and
$\varphi \in C^\infty_c(\Om;\R^m)$.
The rank one convexity of $\xi \mapsto W(x,\xi)$ on $\msd$ and 
the growth assumption \eqref{scisandb2} imply (see, e.g., \cite{D}) that  there exist a positive constant $a_2^W>0$ and a
nonnegative function $b_2^W\in L^{p'}\!(\Om)$, with
$p':=p/(p-1)$, such that
\begin{equation}
\label{estgradxi}
|\partial_\xi W(x,\xi)|\le a_2^W |\xi|^{p-1} + b_2^W(x)\hbox{  for all
 }(x,\xi)\in \Om{\times}  \msd,
\end{equation}
where $\partial_\xi W\colon \Om {\times} \msd \to \msd$ denotes
the partial gradient of $W$ with respect to
$\xi$.
By \eqref{scisandb2} and  \eqref{estgradxi} the functional
$\ws$, defined for all $\Phi \in L^p(\Om; \msd)$ by
$$
\ws(\Phi):=\int_\Om W(x,\Phi(x)) \, dx,
$$
is of class $C^1$ on $L^p(\Om;\msd)$, and its differential
$\partial \ws:L^p(\Om;\msd) \to L^{p'}\!(\Om;\msd)$ is given by
\begin{equation*}
\langle \partial \ws(\Phi),\Psi \rangle =
\int_\Om \partial_\xi W(x,\Phi(x)) \Psi(x)\, dx
\qquad\hbox{for every }
\Phi,\Psi \in L^p(\Om;\msd),
\end{equation*}
where $\langle \cdot,\cdot \rangle$ denotes the duality
pairing between the spaces
$L^{p'}\!(\Om;\msd)$ and $L^p(\Om;\msd)$.

\par\medskip\noindent
{\bf The body forces.} Since we consider only conservative body and surface forces, it is convenient to describe them by means of their potentials, that will be denoted by $F$ and $G$, respectively. Let $q>1$ be fixed.
The density of the applied
body forces per unit volume in the
reference configuration relative to the deformation
$u$ at time $t \in [0,T]$
is then given by $\partial_z F(t,x,u(x))$. We assume that
$F\colon[0,T]{\times} \Om {\times} \R^m \to \R$ satisfies the following conditions:
\begin{eqnarray*}
&x\mapsto F(t,x,z)
\text{ is } \leb^n \text{ measurable on }
\Om\text{ for every } (t,z)\in [0,T]{\times} \R^m, \\
&z \mapsto F(t,x,z) \text{ belongs to } C^1(\R^m)
\text{ for every } (t,x)\in [0,T] {\times} \Om.
\end{eqnarray*}

We assume that for every $t\in[0,T]$ the functional
\begin{equation}\label{bodyenergy}
\fs(t)(u):= \int_\Om F(t,x,u(x))\,dx
\end{equation}
is of class $C^1$ on $L^q(\Om;\R^m)$, with differential 
$\partial \fs(t)\colon L^q(\Om;\R^m)\to L^{q'}\!(\Om;\R^m)$,
$q':=\frac q{q-1}$,
given by
$$
\langle \partial \fs(t)(u),v \rangle =
\int_\Om \partial_z F(t,x,u(x)) v(x)\, dx
\qquad\hbox{for every }
u,v \in L^q(\Om;\R^m),
$$
where $\langle \cdot,\cdot \rangle$ denotes now the duality
pairing between
$L^{q'}\!(\Om;\R^m)$ and $L^q(\Om;\R^m)$.
We assume also the following semicontinuity condition:
$$
\fs(t)(u)\ge\limsup_{k\to \infty} \fs(t)(u_k)
$$
for every $u_k$, $u\in L^q(\Om;\R^m)$ such that $u_k\to u$ a.e.\ on $\Om$.

As for the regularity with respect to time, we assume that there exist a constant 
$\dot q<q$ and, for a.e.\ $t\in [0,T]$, a functional 
$\dot\fs(t)\colon  L^{\dot q}(\Om;\R^m)\to\R$ of class $C^1$, with differential 
$\partial \dot\fs(t)\colon L^{\dot q}(\Om;\R^m)\to L^{\dot q'}\!(\Om;\R^m)$, $\dot{q}':=\frac{\dot{q}}{\dot{q}-1}$, such that for every $u$, $v\in L^q(\Om;\R^m)$
the functions $t\mapsto \dot\fs(t)(u)$ and
$t\mapsto \langle\partial\dot\fs(t)(u),v\rangle$ are integrable on $[0,T]$, and
\begin{eqnarray*}
&
\displaystyle
\fs(t)(u) = \fs(0)(u)+ \int_0^t \dot\fs(s)(u) \,ds,\\
&
\displaystyle
\langle\partial\fs(t)(u),v\rangle= \langle\partial\fs(0)(u),v\rangle+
\int_0^t \langle\partial\dot\fs(s)(u),v\rangle \,ds \quad\quad
\end{eqnarray*}
for every $t\in[0,T]$.

We assume that $\fs(t)$, $\partial \fs(t)$, $\dot\fs(t)$, and $\partial\dot\fs(t)$ satisfy the following growth conditions for every $u$, $v\in L^q(\Om;\R^m)$ and
a.e.\ $t\in[0,T]$:
\begin{eqnarray*}
&a_0^\fs\|u\|_q^q - b_0^\fs \le -\fs(t)(u)\le a_1^\fs\|u\|_q^q + b_1^\fs, 
\\
&|\langle\partial\fs(t)(u),v\rangle| \le \big(a_2^\fs \|u\|_q^{q-1} + b_2^\fs\big)\|v\|_q,
\\
& |\dot\fs(t)(u)|\le a_3^\fs(t)\|u\|_{\dot{q}}^{\dot{q}} + b_3^\fs(t), 
\\
& |\langle\partial\dot\fs(t)(u),v\rangle|\le \big(a_4^\fs(t)\|u\|_{\dot{q}}^{\dot{q}-1}
+ b_4^\fs(t)\big)\|u\|_{\dot{q}} ,
\end{eqnarray*}
where  $a_0^{\fs}>0$, $a_1^{\fs}>0$, $a_2^{\fs}>0$,
$b_0^{\fs} \ge 0$, $b_1^{\fs} \ge 0$, and
$b_2^{\fs} \ge 0$ are constants,
$a_3^{\fs}$, $a_4^{\fs}$, $b_3^{\fs}$, and
$b_4^{\fs}$
are nonnegative functions in $L^1([0,T])$, and $\|\cdot\|_s$ denotes the norm in $L^s$.

\par
The positivity of $a_0^\fs$ is a crucial assumption for the coerciveness of the elastic
energy defined in \eqref{elener} below. Since it implies that $-\fs(t)(u)$ is large
for large values of $\|u\|_q$, the forces cannot send  portions
of the body to infinity, even if the cracks happen to split it into several pieces. This allows to obtain an existence result of the quasistatic evolution
for arbitrarily large times.

\par\medskip\noindent
{\bf The surface forces.}
The density of the surface forces on $\partial_S \Om$ at
time $t$ under the
deformation $u$ is given by $\partial_z G(t,x,u(x))$,
where $G\colon[0,T] {\times} \partial_S \Om {\times} \R^m \to \R$
is such that
\begin{eqnarray*}
&x \mapsto G(t,x,z)
\text{ is }  \hn\text{{-}measurable on } \partial_S\Om \text{ for every } (t,z) \in  [0,T] {\times}\R^m,
\\
&
z \mapsto G(t,x,z) \text{ belongs to } C^1(\R^m)\text{ for every } (t,x) \in [0,T] {\times} \partial_S\Om.
\end{eqnarray*}

Let us fix an exponent $r$, whose value is related to the trace operators on
Sobolev spaces: if $p<n$ we suppose
that $p \le r \le \frac{p}{n-p}$, while if $p \ge n$
we suppose only $p \le r$. 
We assume that for every $t\in[0,T]$ the functional
\begin{equation}\label{tracener}
\gs(t)(u):= \int_{\partial_S\Om} G(t,x,u(x))\,d\hn(x)
\end{equation}
is of class $C^1$ on $L^r(\partial_S\Om;\R^m)$, with differential 
$\partial \gs(t)\colon L^r(\partial_S\Om;\R^m)\to L^{r'}\!(\partial_S\Om;\R^m)$, 
$r':=\frac r{r-1}$, given by
$$
\langle \partial \gs(t)(u),v \rangle =
\int_\Om \partial_z G(t,x,u(x)) v(x)\, d\hn(x),
\qquad\hbox{for every }
u,v \in L^r(\partial_S\Om;\R^m),
$$
where $\langle \cdot,\cdot \rangle$ denotes now the duality
pairing between
$L^{r'}\!(\Om;\R^m)$ and $L^r(\Om;\R^m)$.

As for the regularity with respect to time, we assume that for a.e.\ $t\in [0,T]$
there exists a functional 
$\dot\gs(t)\colon  L^r(\partial_S\Om;\R^m)\to\R$ of class $C^1$, with differential 
$\partial \dot\gs(t)\colon L^r(\partial_S\Om;\R^m)\to L^{r'}\!(\partial_S\Om;\R^m)$,
such that for every $u$, $v\in L^r(\partial_S\Om;\R^m)$
the functions $t\mapsto \dot\gs(t)(u)$ and
$t\mapsto \langle\partial\dot\gs(t)(u),v\rangle$ are integrable on $[0,T]$, and
\begin{eqnarray*}
&
\displaystyle
\gs(t)(u) = \gs(0)(u)+ \int_0^t \dot\gs(s)(u) \,ds,\\
&
\displaystyle
\langle\partial\gs(t)(u),v\rangle= \langle\partial\gs(0)(u),v\rangle+
\int_0^t \langle\partial\dot\gs(s)(u),v\rangle \,ds \quad\quad
\end{eqnarray*}
for every $t\in[0,T]$.

We assume that $\gs(t)$, $\partial \gs(t)$, $\dot\gs(t)$, and $\partial\dot\gs(t)$ satisfy the following growth conditions for every $u$, $v\in L^r(\partial_S\Om;\R^m)$ and
a.e.\ $t\in[0,T]$:
\begin{eqnarray*}
&-a_0^\gs\|u\|_{r,\partial_S\Om} - b_0^\gs \le -\gs(t)(u)\le a_1^\gs\|u\|_{r,\partial_S\Om}^r+ b_1^\gs, 
\\
&|\langle\partial\gs(t)(u),v\rangle| \le \big(a_2^\gs \|u\|_{r,\partial_S\Om}^{r-1} + b_2^\gs\big)\|v\|_r,
\\
& |\dot\gs(t)(u)|\le a_3^\gs(t)\|u\|_{r,\partial_S\Om}^r + b_3^\gs(t), 
\\
& |\langle\partial\dot\gs(t)(u),v\rangle|\le \big(a_4^\gs(t)\|u\|_{r,\partial_S\Om}^{r-1}
+ b_4^\gs(t)\big)\|v\|_{r,\partial_S\Om},
\end{eqnarray*}
where  $a_0^{\gs}$, $a_1^{\gs}$, $a_2^{\gs}$,
$b_0^{\gs}$, $b_1^{\gs}$, and
$b_2^{\gs}$ are nonnegative constants,
$a_3^{\gs}$, $a_4^{\gs}$, $b_3^{\gs}$, and
$b_4^{\gs}$
are nonnegative functions in $ L^1([0,T])$, and $\|\cdot\|_{r,\partial_S\Om}$ denotes the norm in
$L^r(\partial_S\Om;\R^m)$.

\par\medskip\noindent
{\bf Configurations with finite energy.}
The deformations on the boundary $\partial_D \Om$ are given by
(the traces of) functions
$\psi \in W^{1,p}(\Om;\R^m)\cap L^q(\Om;\R^m)$. Given a crack $\Gamma \in \radm$ and 
a boundary deformation $\psi$,
the set of {\it admissible deformations with finite energy} 
relative to $(\psi,\Gamma)$
is defined by
\begin{equation*}
AD(\psi,\Gamma):= \{u \in GSBV^p_q(\Om;\R^m):
S(u) \tsub \Gamma,\ u=\psi\,\ \hn\text{-a.e.\ on }\partial_D\Om\setminus\Gamma\}.
\end{equation*}

\par
Note that, if $u \in GSBV^p_q(\Om;\R^m)$, then
$\ws(\nabla u)<+\infty$ and
$|\fs(t)(u)|<+\infty$ for all $t \in [0,T]$. Moreover,
since $\Gamma \in \radm$, $S(u) \tsub \Gamma\tsub\overline\Om_B$, and 
$\partial_S\Om\cap \overline\Om_B=\emptyset$,
we have that
$\gs(t)(u)$ is well defined and
$|\gs(t)(u)|<+\infty$ for all $t \in [0,T]$ (see \cite[Section~3]{DMFT} for the details).
Notice that  there always exists a deformation without
crack which satisfies the boundary condition,
namely the function $\psi$ itself. This means that, if some cracks appear, this is because they are energetically convenient with respect to the elastic solution, and not because the presence of a crack is the only way to match the boundary condition.

\par\medskip\noindent
{\bf The total energy.}
For every $t\in [0,T]$ the total energy of the
configuration $(u,\Gamma)$, with $u \in AD(\psi,\Gamma)$,
is given by
\begin{equation}
\label{totalener}
\Es(t)(u,\Gamma):=\Eb(t)(u)+\Esup(\Gamma),
\end{equation}
where the {\it surface energy\/} $\Esup$ is defined in  \eqref{crackener}, while the {\it elastic energy\/}
$\Eb(t)$ is given by
\begin{equation}
\label{elener}
\Eb(t)(u):=\ws(\nabla u)-\fs(t)(u)-\gs(t)(u),
\end{equation}
with $\ws$, $\fs(t)$, and $\gs(t)$ defined in
\eqref{bws}, \eqref{bodyenergy},
and \eqref{tracener}, respectively.
\par\medskip\noindent
{\bf The time dependent boundary deformations.}
We will consider boundary deformations $\psi(t)$ such that
$$
t \mapsto \psi(t) \in \,AC([0,T];
W^{1,p}(\Om;\R^m) \cap L^q(\Om;\R^m)),
$$
so that
\begin{eqnarray*}
&t \mapsto \dot\psi(t) \in L^1([0,T];
W^{1,p}(\Om;\R^m) \cap L^q(\Om;\R^m)),
\\
&
t \mapsto \nabla\dot\psi(t)\in L^1([0,T];L^p(\Om;\msd)).
\end{eqnarray*}

\par\medskip\noindent
{\bf Quasistatic evolution.}
The notion of quasistatic evolution of a cracked configuration is made precise by the
following definition.

\begin{definition}
\label{qsedef}
A {\it quasistatic evolution\/} with boundary deformation $t \mapsto \psi(t)$ is a
function $t \mapsto (u(t),\Gamma(t))$ from
$[0,T]$ to $GSBV^p_q(\Om;\R^m){\times} \radm$
with the following properties:
\begin{itemize}
\vskip4pt
\item[(a)]
{\it global stability}:\,\,for all $t\in [0,T]$ we have
$u(t) \in AD(\psi(t),\Gamma(t))$  and
$$
\Es(t)(u(t),\Gamma(t)) = \min \{\Es(t)(v,\Gamma):\Gamma\in \radm,\ \Gamma(t) \tsub \Gamma,\ v \in
AD(\psi(t),\Gamma)\};
$$
\vskip4pt
\item[(b)]
{\it irreversibility}:\,\,
$\Gamma(s)\tsub\Gamma(t)$ whenever
$0\le s<t \le T$;
\vskip4pt
\item[(c)]
{\it energy balance}:\,\,
the function $t \mapsto E(t):=\Es(t)(u(t),\Gamma(t))$ is
absolutely continuous on  $[0,T]$ and
\begin{equation}\label{nondissqse}
\begin{array}{c}
\dot{E}(t) = \langle\partial\ws(\nabla u(t)),
\nabla \dot\psi(t)\rangle
-\langle\partial \fs(t)(u(t)),\dot\psi(t)\rangle-
\dot{\fs}(t)(u(t)) \vspace{2pt}\\
-\langle\partial \gs(t)(u(t)),\dot\psi(t)\rangle-
\dot{\gs}(t)(u(t))
\end{array}
\end{equation}
for a.e.  $t\in [0,T]$.
\end{itemize}
\end{definition}

If the solution is sufficiently regular, an integration by parts shows that the right-hand side
of \eqref{nondissqse} represents the power of all external forces acting on the body,
including the unknown forces on $\partial_D\Om\!\setminus\!\Gamma(t)$ that produce 
the imposed boundary deformation $\psi(t)$
(see \cite[Section~3.9]{DMFT} for details).

\par\medskip\noindent
{\bf The existence result.}
We are interested in quasistatic evolutions with a prescribed initial condition
$(u_0,\Gamma_0)$ with $\Gamma_0\in \radm$ and $u_0\in AD(\psi(0),\Gamma_0)$. 
Since in the definition we require that the global stability condition
is satisfied for every time, a necessary condition for the solvability of the initial
 value problem is that
\begin{equation}\label{mininit}
\Es(0)(u_0,\Gamma_0)\le \Es(0)(u,\Gamma)
\end{equation}
for every $ \Gamma\in \radm$, with
$\Gamma_0\tsub\Gamma$, and every $u\in AD(\psi(0),\Gamma)$.

The next Theorem, proved in \cite{DMFT}, establishes
the existence of a quasistatic evolution with prescribed initial and boundary conditions.

\begin{theorem}
\label{qse}
Let $\Gamma_0\in \radm$ and $u_0\in AD(\psi(0),\Gamma_0)$. Assume that 
\eqref{mininit} is satisfied.
Then there exists a quasistatic evolution with boundary deformation $\psi(t)$ such that
$(u(0),\Gamma(0))=(u_0,\Gamma_0)$.
\end{theorem}

\section{Qualitative properties of the quasistatic crack growth}
\label{qualsec}
In this section we consider a quasistatic crack growth $t \mapsto (u(t),\Gamma(t))$
with boundary deformation $t \mapsto \psi(t)$ according to Definition~\ref{qsedef}.
In the sequel we will derive some qualitative properties of the cracks $\Gamma(t)$
and of the deformations $u(t)$. 

\subsection{Left continuous envelope of the crack}
For every $t \in {]0,T]}$ we define $\Gamma^-(t)$ as the rectifiable set 
$$
\Gamma^-(t):= \bigcup_{k}\Gamma(s_k),
$$
where $(s_k)$ is a sequence converging to $t$ with $s_k < t$ for every $k$. 
We define $\Gamma^-(0):=\Gamma(0)$. Since $\Gamma(\cdot)$ is increasing in time,
it turns out that $\Gamma^-(t)$ is independent of the choice of the sequence $(s_k)$.

\begin{proposition}\label{leftcont}
For all $t \in [0,T]$ we can find $v(t) \in AD(\psi(t),\Gamma^-(t))$ such that
$t \mapsto (v(t),\Gamma^-(t))$ is a quasistatic evolution  with boundary deformation~$\psi(t)$.
\end{proposition}

\begin{dimo}
If follows from the definition that $\Gamma^-(t)$ is increasing with respect to~$t$.
Since $t \mapsto \Gamma(t)$ is increasing, it turns out that $\Gamma^-(t)=\Gamma(t)$
for all $t \in [0,T]$ except for a countable set $C\subset {]0,T]}$. For every
$t \in [0,T]\!\setminus\!C$ we set
$v(t):=u(t)$,
while for $t \in C$ we take as $v(t)$ any minimizer of
$$
\min\{ \Eb(t)(u)\,:\, u \in AD(\psi(t),\Gamma^-(t))\}.
$$

\par
Let us check that the three conditions for quasistatic evolutions hold for the function 
$t \mapsto (v(t),\Gamma^-(t))$. Irreversibility has already been proved. 
As for global stability, let $\Gamma$ be such that $\Gamma^-(t) \tsub \Gamma$, 
and let $v \in AD(\psi(t),\Gamma)$. We notice that $\Gamma(s) \tsub \Gamma$ 
for every $s<t$, 
so that by the global stability of $(u(s),\Gamma(s))$ we get
$$
\Eb(s)(u(s))+\Esup(\Gamma(s))=
\Es(s)(u(s),\Gamma(s))
\le \Es(s)(v-\psi(t)+\psi(s),\Gamma).
$$ 
Let $(s_k)$ with $s_k \to t$ and $s_k<t$ for every $k$: up to a subsequence we 
have that $u(s_k) \weak \tilde u$ weakly in $GSBV^p_q(\Om;\R^m)$ for some 
$\tilde u \in AD(\psi(t),\Gamma^-(t))$. By lower semicontinuity the previous 
inequality gives
$$
\Eb(t)(\tilde u)+\Esup(\Gamma^-(t)) \le \Es(t)(v,\Gamma),
$$
so that by the minimality of $v(t)$
$$
\Es(t)(v(t),\Gamma^-(t))=
\Eb(t)(v(t))+\Esup(\Gamma^-(t)) \le \Es(t)(v,\Gamma).
$$

\par
Let us come to the energy balance. Since $(v(t),\Gamma^-(t))=(u(t),\Gamma(t))$ 
for all $t$ up to a countable set, it is sufficient to prove that 
$\Es(t)(v(t),\Gamma^-(t))=\Es(t)(u(t),\Gamma(t))$ for all $t \in [0,T]$.
By global minimality for every $s<t$ we have
$$
\Es(s)(u(s),\Gamma(s)) \le \Es(s)(v(t)-\psi(t)+\psi(s),\Gamma^-(t))
$$
and so, since $t\mapsto\Es(t)(u(t),\Gamma(t))$ is continuous by definition, we get
for $s \to t$
$$
\Es(t)(u(t),\Gamma(t)) \le \Es(t)(v(t),\Gamma^-(t)).
$$
The opposite inequality comes from the global stability of $(v(t),\Gamma^-(t))$, so that the proposition is proved.
\end{dimo}

We can also consider the right continuous envelope of the crack.
For every $t \in {[0,T[}$ we define $\Gamma^+(t)$ as the rectifiable set 
$$
\Gamma^+(t):= \bigcap_{k}\Gamma(s_k),
$$
where $(s_k)$ is a sequence converging to $t$ with $s_k > t$ for every $k$. 
We define $\Gamma^+(T):=\Gamma(T)$. Since $\Gamma(\cdot)$ is increasing in time,
it turns out that $\Gamma^+(t)$ is independent of the choice of the sequence $(s_k)$.

\begin{proposition}
For all $t \in [0,T]$ we can find $v(t) \in AD(\psi(t),\Gamma^+(t))$ such that
$t \mapsto (v(t),\Gamma^+(t))$ is a quasistatic evolution with boundary deformation
$\psi(t)$.
\end{proposition}

\begin{dimo}
As in the proof of Proposition~\ref{leftcont}
we find that $\Gamma^+(t)=\Gamma(t)$
for all $t \in [0,T]$ except for a countable set $C\subset {[0,T[}$, and
for every $t \in C$ we take as $v(t)$ any minimizer of
$$
\min\{ \Eb(t)(u)\,:\, u \in AD(\psi(t),\Gamma^+(t))\}.
$$

\par
Let us check global stability and energy balance for the function 
$t \mapsto (v(t),\Gamma^+(t))$.
As for global stability, let $\Gamma$ be such that $\Gamma^+(t) \tsub \Gamma$, 
and let $v \in AD(\psi(t),\Gamma)$. For every $s>t$
the global stability of $(u(s),\Gamma(s))$ gives
$$
\Es(s)(u(s),\Gamma(s))
\le \Es(s)(v-\psi(t)+\psi(s),\Gamma\cup\Gamma(s)).
$$ 
Let $(s_k)$ with $s_k \to t$ and $s_k>t$ for every $k$: up to a subsequence 
we have that $u(s_k) \weak \tilde u$ weakly in $GSBV^p_q(\Om;\R^m)$ for 
some $\tilde u \in AD(\psi(t),\Gamma^+(t))$. By lower semicontinuity the previous inequality gives
$$
\Eb(t)(\tilde u)+\Esup(\Gamma^+(t)) \le \Es(t)(v,\Gamma),
$$
so that by the minimality of $v(t)$
$$
\Es(t)(v(t),\Gamma^+(t))=\Eb(t)(v(t))+\Esup(\Gamma^+(t)) \le \Es(t)(v,\Gamma).
$$

\par
Let us come to the energy balance. As in Proposition~\ref{leftcont} it is sufficient to prove that 
$\Es(t)(v(t),\Gamma(t))=\Es(t)(v(t),\Gamma^+(t))$ for all $t \in [0,T]$. 
The global stability of $(v(t),\Gamma^+(t))$ gives, for every $s>t$,
$$
\Es(t)(v(t),\Gamma^+(t)) \le \Es(s)(u(s)-\psi(s)+\psi(t),\Gamma(s))
\le \Es(s)(u(s),\Gamma(s))+\omega(s),
$$
where $\omega(s)\to 0$ as $s\to t$ (see 
the estimates of \cite[Section~5]{DMFT}). By the continuity of 
 $t\mapsto\Es(t)(u(t),\Gamma(t))$ we obtain
for $s \to t$
$$
\Es(t)(v(t),\Gamma^+(t)) \le\Es(t)(u(t),\Gamma(t)).
$$
The opposite inequality comes from the global stability of $(u(t),\Gamma(t))$, 
so that the proposition is proved.
\end{dimo}

\subsection{Structure of the crack set}
Francfort and Larsen, in their approach to quasistatic crack growth \cite{FL}, define
the crack set at time $t$ as the union of the jump sets of the displacements at previous times. 
In the next proposition we prove a similar structure result for any quasistatic evolution 
according to Definition~\ref{qsedef}.

\begin{proposition}
\label{strctgammaprop}
There exists a countable and dense set $D \subset [0,T]$ such that for all $t \in [0,T]$
\begin{equation}
\label{structureGamma}
\Gamma(t)\teq\Gamma(0)\cup \bigcup_{s \in D,\,s \le t } \Sg{\psi(s)}{u(s)},
\end{equation}
where $\teq$ means equality up to a set of $\hn$-measure zero and
\begin{equation}
\label{Spsi}
S^\psi(u):=S(u)\cup\{x\in\partial_D\Om:u(x)\neq\psi(x)\}
\end{equation}
for every $\psi\in W^{1,p}(\Om;\R^m)$ and $u\in GSBV^p_q(\Om;\R^m)$.
\end{proposition}

\begin{dimo}
By \cite[Lemma~4.12]{DMFT} there exists a sequence of subdivisions
$(t^i_k)_{0 \le i \le i_k}$ of $[0,T]$, with
\begin{equation}
\label{subdiv}
0=t^0_k < t^1_k < \dots < t^{i_k-1}_k < t^{i_k}_k=T
\quad\text{and}\quad
\lim_{k \to \infty} \max_{1 \le i \le i_k} (t^i_k-t^{i-1}_k)=0,
\end{equation}
 such that
\begin{equation}
\label{thetaerror}
\lim_{k \to \infty} \sum_{i=1}^{i_k}
\left|(t^i_k-t^{i-1}_k)\dot E(t^i_k)- \int_{t^{i-1}_k}^{t^i_k} \dot E(t)\,dt \right|=0.
\end{equation}
 Moreover, we can assume that $\{t^i_k: 0 \le i \le i_k\} \subset 
 \{t^i_{k+1}:0 \le i \le i_{k+1}\}$ for all $k$. Let us define
$$
I_{\infty}:= \{t^i_k: k\in\N,\ 0 \le i \le i_k\}.
$$
For all $t \in [0,T]$ we set
$$
\tilde{\Gamma}(t):= \Gamma(0)\cup \bigcup_{s \in I_{\infty},\, s \le t} \Sg{\psi(s)}{u(s)},
$$
and
$$
\tilde{E}(t):=\Eb(t)(u(t))+\Esup(\tilde{\Gamma}(t)).
$$
Notice that $\tilde{\Gamma}(0)=\Gamma(0)$ and $\tilde{\Gamma}(t) \tsub \Gamma(t)$ 
for all $t \in [0,T]$. As a consequence $\tilde{E}(0)=E(0)$ and $\tilde{E}(t) \le E(t)$ for all $t \in [0,T]$. 
\par
Let us fix $t \in I_{\infty}$. For every $k$ sufficiently large there exists
$j_k \in \{0,\dots,i_k\}$ such that
$t=t^{j_k}_k$. 
For every $i=1,\dots,j_k-1$, by the global stability of $(u(t^i_k),\Gamma(t^i_k))$ we get
\begin{eqnarray*}
\Eb(t^i_k)(u(t^i_k))\!\!&\le&\!\!
\Eb(t^i_k)(u(t^{i+1}_k)-\psi(t^{i+1}_k)+\psi(t^i_k))+
\Esup(\Sg{\psi(t^{i+1}_k)}{u(t^{i+1}_k)} \setminus \Gamma(t^i_k)) \\
\!\!&\le&\!\!
\Eb(t^i_k)(u(t^{i+1}_k)-\psi(t^{i+1}_k)+\psi(t^i_k))+
\Esup(\Sg{\psi(t^{i+1}_k)}{u(t^{i+1}_k)} \setminus \tilde \Gamma(t^i_k)),
\end{eqnarray*}
so that
$$
\Eb(t^i_k)(u(t^i_k))+\Esup(\tilde \Gamma(t^i_k)) \le
\Eb(t^i_k)(u(t^{i+1}_k)-\psi(t^{i+1}_k)+\psi(t^i_k))+
\Esup(\tilde \Gamma(t^{i+1}_k)).
 $$
Using the error estimates contained in formula (7.46) in the proof of  
\cite[Theorem~3.15]{DMFT},
we deduce that there exists a sequence $e_k(t) \to 0$ such that
\begin{equation*}
\tilde{E}(t) \ge E(0)+
\sum_{j=1}^{j_k} (t^i_k-t^{i-1}_k) \dot E(t^i_k)-e_k(t).
\end{equation*}
By \eqref{thetaerror} we deduce that for all
$t \in I_{\infty}$
$$
\tilde{E}(t) \ge E(0)+\int_0^t  \dot E(s) \,ds=E(t).
$$
Since we have already proved the opposite inequality, we obtain $\tilde{E}(t)=E(t)$ 
for all $t \in I_\infty$. As a consequence  for all $t \in I_\infty$ we have 
$\Gamma(t) \teq \tilde{\Gamma}(t)$, hence
\begin{equation}
\label{representation}
\Gamma(t)\teq \Gamma(0)\cup \bigcup_{s \in I_{\infty},\, s \le t} \Sg{\psi(s)}{u(s)}.
\end{equation}
Equality \eqref{representation} extends to all continuity
points of $t\mapsto\hn(\Gamma(t))$.
\par
Let us set $D:= J\cup I_\infty $, where $J$ is the (at most countable) set of jumps point of 
$t\mapsto\hn(\Gamma(t))$.
Notice that 
\begin{equation}
\label{gammastructure}
\Gamma(t)\teq\Gamma(0)\cup \bigcup_{s \in D,\, s \le t} \Sg{\psi(s)}{u(s)}
\end{equation}
for every $t\in[0,T]\!\setminus\! J$.  To conclude the proof it is thus sufficient to show that
\eqref{gammastructure} holds also for $t \in J$. Let $\hat \Gamma(t)$ be the right-hand 
side of \eqref{gammastructure}, and let us prove that $\hat \Gamma(t) \teq \Gamma(t)$.
Let $t\in J$ and let $t_k \in I_{\infty}$ with $t_k \to t$ and $t_k\le t$ for every $k$.
Since $v_k:=u(t)-\psi(t)+\psi(t_k)$ satisfies 
$\Sg{\psi(t_k)}{v_k}=\Sg{\psi(t)}{u(t)} \subset \hat{\Gamma}(t)$ and $\Gamma(t_k) \tsub \hat{\Gamma}(t)$,
by the global stability of $(u(t_k),\Gamma(t_k))$ we get
$$
E(t_k)=\Eb(t_k)(u(t_k))+\Esup(\Gamma(t_k)) \le \Eb(t_k)(v_k)+\Esup(\hat{\Gamma}(t)).
$$
Letting $k \to \infty$, by the continuity of $t\mapsto E(t)$ we obtain
$$
E(t) \le \Eb(t)(u(t))+\Esup(\hat{\Gamma}(t)),
$$
which implies $\Esup(\Gamma(t)) \le \Esup(\hat{\Gamma}(t))$.
Since $\hat{\Gamma}(t) \tsub \Gamma(t)$, we conclude that $\Gamma(t) \teq \hat{\Gamma}(t)$.
\end{dimo}

The structure result can be improved under suitable convexity assumptions  for the bulk 
energy and for the potentials of the applied forces.

\begin{lemma}
\label{strctgammaconvex}
 For every $t \in [0,T]$ assume that $\ws$ and $-\gs(t)$ are convex and that 
 $-\fs(t)$ is strictly convex.
Then for every countable dense set $D \subset [0,T]$ containing the jump points of 
$t\mapsto\hn(\Gamma(t))$ we have
\begin{equation}
\label{structureGammaconvex}
\Gamma(t)\teq\Gamma(0)\cup  \bigcup_{s \in D,\, s \le t} \Sg{\psi(s)}{u(s)}.
\end{equation}
\end{lemma}

\begin{dimo}
By Proposition \ref{strctgammaprop} we know that there exists a countable and
dense set $D' \subset [0,T]$, containing the set of jumps of 
$t\mapsto\hn(\Gamma(t))$, such that
\begin{equation}
\label{structuregamma}
\Gamma(t)\teq\Gamma(0)\cup  \bigcup_{s \in D',\, s \le t} \Sg{\psi(s)}{u(s)}.
\end{equation}
Let $D$ be a countable and dense subset of $[0,T]$, containing the set of jumps of 
$t\mapsto\hn(\Gamma(t))$, and let
$$
\tilde{\Gamma}(t):=\Gamma(0)\cup \bigcup_{s \in D,\, s \le t} \Sg{\psi(s)}{u(s)}.
$$
Recall that $\tilde{\Gamma}(t) \tsub \Gamma(t)$ for all $t \in [0,T]$ because 
$u(t) \in AD(\psi(t),\Gamma(t))$ for all $t \in [0,T]$. 
Let us show that the opposite inclusion
holds. By \eqref{structuregamma}, it is sufficient to prove that for all $s \in D'$ with
$0<s \le t$ we have $\Sg{\psi(s)}{u(s)} \tsub \tilde{\Gamma}(t)$. If $s$ is a jump point
of $t\mapsto\hn(\Gamma(t))$, we have $s \in D$ and the inclusion follows. If $s$ is not a jump point
of $t\mapsto\hn(\Gamma(t))$, let us consider $s_k \in D$ with $s_k \to s$ and $s_k\le s$.
We have that (up to subsequences) 
$$
u(s_k) \weak \tilde u 
\qquad
\text{weakly in }GSBV^p_q(\Om;\R^m),
$$
and $\Sg{\psi(s)}{\tilde u} \tsub \tilde \Gamma(s)$, since $\Sg{\psi(s_k)}{u(s_k)} \tsub \tilde \Gamma(s)$ 
for all $k$ (see \cite[Remark~2.9]{DMFT}). Moreover $\tilde u$ is a minimizer of
\begin{equation}
\label{minpbebs}
\min \{ \Eb(s)(v) \,:\, v \in GSBV^p_q(\Om;\R^m),\; \Sg{\psi(s)}{v} \tsub \Gamma(s) \}.
\end{equation}
Indeed for all $v \in GSBV^p_q(\Om;\R^m)$ with $\Sg{\psi(s)}{v} \tsub \Gamma(s)$ we have that
$\Sg{\psi(s_k)}{v-\psi(s)+\psi(s_k)}=\Sg{\psi(s)}{v} \tsub \Gamma(s)$, so that
$$
\Eb(s_k)(u(s_k)) \le \Eb(s_k)(v-\psi(s)+\psi(s_k))+ \Esup \left( \Gamma(s) \setminus \Gamma(s_k) \right).
$$
Taking the limit for $k \to \infty$ we obtain
$$
\Eb(s)(\tilde u) \le \Eb(s)(v),
$$
hence $\tilde u$ is a solution of problem \eqref{minpbebs}.
By the convexity assumption on the bulk energy, 
the solution of problem \eqref{minpbebs} is unique. By the global stability condition for $(u(s),\Gamma(s))$
the function $u(s)$ is a minimizer of \eqref{minpbebs}. This implies $\tilde u=u(s)$, 
which gives $\Sg{\psi(s)}{u(s)} \tsub \tilde{\Gamma}(s) \tsub \tilde{\Gamma}(t)$ and concludes the proof.
\end{dimo}

\subsection{Measurability properties of the deformation}
In the quasistatic evolution $t \mapsto (u(t),\Gamma(t))$, 
the deformation $u(t)$ is, in general, not uniquely determined by $\Gamma(t)$. 
Indeed, if $v(t)$ is another minimizer of $\Eb(t)$ in  $AD(\psi(t),\Gamma(t))$, 
then the global stability condition still holds and the value of $E(t)$ does not change. 
Therefore $t\mapsto (v(t),\Gamma(t))$ is a quasistatic evolution provided \eqref{nondissqse} 
is satisfied by $v(t)$. The following result shows that we can select $t\mapsto v(t)$ 
so that certain measurability properties hold. Similar results for a different evolution problem have been obtained independently in~\cite{Mai}.

\begin{theorem}
\label{meausut}
We can choose $t\mapsto v(t)$ in such a way that $t\mapsto (v(t),\Gamma(t))$ 
is a quasistatic evolution with boundary condition $t\mapsto\psi(t)$ and  
$t \mapsto (\nabla v(t), v(t))$ 
is measurable from $[0,T]$ to $L^p(\Om; \msd) {\times} L^q(\Om;\R^m)$.
\end{theorem}

\begin{dimo}
Let $t^i_k$ and $I_\infty$ be as in Proposition~\ref{strctgammaprop}.
For all $k$ we set
\begin{equation}
\label{vkt}
u_k(t):=u(t^{i-1}_k),
\qquad
\Gamma_k(t):=\Gamma(t^{i-1}_k),
\qquad
\vartheta_k(t):=\dot E(t^{i-1}_k)
\qquad
\text{for }t \in {[t^{i-1}_k,t^i_k[}.
\end{equation}
Arguing as in the proof of \cite[Lemma~6.1]{DMFT}, we deduce from \eqref{thetaerror}
that for all $t \in [0,T]$ there exists $e_k(t) \to 0$ as $k \to \infty$ such that 
$$
\Es(t)(u_k(t),\Gamma_k(t))=
\Es(0)(u_k(0),\Gamma_k(0))+\int_0^t \vartheta_k(s)\,ds
+ e_k(t).
$$
Let us define
$$
\tilde \Gamma(t):= \bigcup_{s \in I_\infty,\, s \le t}\Gamma(s).
$$
Using the definition of $\sigma^p$-convergence, introduced in \cite[Section~4.1]{DMFT}, 
and Proposition~\ref{strctgammaprop}
it is easy to prove that for every $t \in [0,T]$ there exists  $\hat \Gamma(t)$ such that
$$
\Gamma_k(t) \stackrel{\sigma^p}{\to} \hat\Gamma(t)\qquad\hbox{and}\qquad
\tilde \Gamma(t)=\Gamma(0)\cup\hat\Gamma(t).
$$
It is clear that $\Gamma(t)=\tilde \Gamma(t)$ for all $t \in I_\infty$. 
As a consequence of the monotonicity with respect to $t$ we obtain that 
$\Gamma(t)=\tilde \Gamma(t)$ for all $t\in[0,T]$ except for a countable set $C$.

\par
Setting
$$
\vartheta(t):=\limsup_{k \to \infty}\vartheta_k(t),
$$
let us consider the sets
\begin{equation}
\label{ast}
\begin{array}{c}
\as(t):=\{(\nabla v, v):  v \in GSBV^p_q(\Om;\R^m),\   u_{k_j}(t) \weak  v 
\text{ weakly in }GSBV^p_q(\Om;\R^m), 
\vspace{2pt}\\
\qquad\vartheta_{k_j}(t) \to \vartheta(t) \text{ for some sequence $k_j\to\infty$}\},
\end{array}
\end{equation}
where $u_k(t)$ is defined in \eqref{vkt}.
Arguing as in \cite[Section~7]{DMFT} we can prove that for every selection 
$t \mapsto (\nabla \tilde v(t), \tilde v(t)) \in \as(t)$
the function $t\mapsto( \tilde v(t),\tilde \Gamma(t))$ is a quasistatic evolution.
Since $\tilde \Gamma(t)=\Gamma(t)$ for $t\not\in C$, the global stability condition gives
$\Es(\tilde v(t), \tilde \Gamma(t))=\Es(u(t),\Gamma(t))$  for $t\not\in C$, and the 
continuity condition in the energy balance ensures that 
$\Es(\tilde v(t), \tilde \Gamma(t))=\Es(u(t),\Gamma(t))$ for every $t\in[0,T]$. Therefore,
if we define $v(t)=\tilde v(t)$ for $t\not\in C$, and $v(t)=u(t)$ for $t\in C$, the function
$t\mapsto(v(t),\Gamma(t))$ is still a quasistatic evolution, and $t\mapsto(\nabla v(t),v(t))$
has the same measurability properties as $t\mapsto(\nabla \tilde v(t),\tilde v(t))$.

\par
In order to conclude the proof, it is enough to show that we can choose $ \tilde v(t)$ so that 
$t \mapsto (\nabla  \tilde v(t),  \tilde v(t))$ is measurable from $[0,T]$ to 
$L^p(\Om;\msd) {\times} L^q(\Om;\R^m)$. 
\par
To this aim, we notice that $(\nabla  v,  v) \in \as(t)$ if and only if there exists a sequence 
$k_j\to\infty$ such that $\vartheta_{k_j}(t) \to \vartheta(t)$,
$\nabla u_{k_j}(t) \weak \nabla  v$ 
weakly in $L^p(\Om;\msd)$,
and
$u_{k_j}(t) \weak  v$ weakly in $L^q(\Om;\R^m)$.
Indeed the convergence in measure of $u_{k_j}(t)$ required in \eqref{1.2} can be obtained 
from the $GSBV$ compactness theorem (see \cite[Theorem~4.36]{AFP}),
since $\hn(S(u_{k_j}(t))$ is uniformly bounded.
Moreover, in view of the coercivity estimates proved in \cite{DMFT}, there exists a 
bounded closed convex set $\ks \subset L^p(\Om;\msd) {\times} L^q(\Om;\R^m)$ such that 
$$
(\nabla u_k(t), u_k(t)) \in \ks
$$
for all $t \in [0,T]$ and for all $k$.
On $\ks$ we consider the weak topology on $L^p(\Om;\msd) {\times} L^q(\Om;\R^m)$, s
o that $\ks$ is a compact metrizable space.
In Lemma \ref{meassel} below we will prove that
\begin{itemize}
\vskip4pt
\item[(a)] $\as(t)$ is closed in the weak topology for every $t \in [0,T]$;
\vskip4pt
\item[(b)] the set $\{t \in [0,T]\,:\, \as(t) \cap U \not= \emptyset\}$ is measurable for every 
open set $U$ in the weak topology of $L^p(\Om;\msd) {\times} L^q(\Om;\R^m)$.
\end{itemize}

\par
Then we can apply the Aumann-von Neumann selection theorem 
(see, e.g., \cite[Theorem III.6]{CV}) and we obtain that we can select 
$(\nabla \tilde v(t), \tilde v(t)) \in \as(t)$ in such a way that 
$t \mapsto (\nabla  \tilde v(t),  \tilde v(t))$ is
measurable from $[0,T]$ to $L^p(\Om;\msd) {\times} L^q(\Om;\R^m)$ 
endowed with the weak topology. The measurability
with respect to the strong topology follows from the Pettis theorem 
(see, e.g.,  \cite[Chapter~V, Section~4]{Yo}) .
\end{dimo}

In the rest of this subsection we prove the lemma concerning the measurability 
of the set valued map $t \mapsto \as(t)$ used in the proof of Theorem \ref{meausut}. 
We settle the problem in the context of compact metric spaces. 
\par
Conditions (a) and (b) in the proof of Theorem  \ref{meausut} follow from the next lemma, 
applied to $X=\ks$ and
$f_k(t)=(\nabla u_k(t), u_k(t))$, with a metric $d$ inducing on $\ks$ the weak topology of 
$L^p(\Om;\msd) {\times} L^q(\Om;\R^m)$.
\begin{lemma}
\label{meassel}
Let $(X,d)$ be a compact metric space, let $f_k\colon[0,T] \to X$ be a sequence of 
measurable functions, and let $\vartheta_k$ and $\vartheta$ be measurable functions from $[0,T]$ to $\R$.
For all $t \in [0,T]$ let
\begin{equation}
\label{ft}
\as(t):= \{x \in X: \text{there exists } k_j \to \infty \text{ such that }f_{k_j}(t)\to x
\text{ and }\vartheta_{k_j}(t)\to \vartheta(t)\}.
\end{equation}
Then
\begin{itemize}
\vskip4pt
\item[$(a)$] $\as(t)$ is closed for all $t \in [0,T]$;
\vskip4pt
\item[$(b)$]  the set $\{t \in [0,T]\,:\, 
\as(t) \cap U \not= \emptyset\}$ is measurable for every open set $U \subset X$.
\end{itemize}
\end{lemma}

\begin{dimo}
Let us fix $t \in [0,T]$ and let us prove that $\as(t)$ is closed
in $X$. Let $x_j \in \as(t)$ with $x_j \to x$ in $X$. Since $x_j \in \as(t)$, we can find
 $k_j \ge j$ such that
$$
d(x_j,f_{k_j}(t)) \le \tfrac{1}{j}
\quad\text{and}\quad
|\vartheta_{k_j}(t)-\vartheta(t)|\le \tfrac{1}{j}.
$$
Clearly
$f_{k_j}(t)\to x$ and
$\vartheta_{k_j}(t)\to \vartheta(t)$,
hence  $x \in \as(t)$ by the very definition of $\as(t)$. 
This proves that $\as(t)$ is closed in $X$.

\par
In view of \cite[Theorem III.9]{CV}, in order to prove $(b)$ it is sufficient to show that
for all $x \in X$ the function $t \mapsto d(x,\as(t))$ is measurable. For every $j$ we define
$$
\as_j(t):=\{x \in X: \text{there exists } k\ge j \text{ such that }d(f_k(t),x)\le\tfrac1j
\text{ and }|\vartheta_k(t)-\vartheta(t)|\le\tfrac1j\},
$$
and we observe that
$
\as(t)=\bigcap_{j}\as_j(t)
$.
We claim that
\begin{equation}
\label{dxft}
d(x,\as(t))=\textstyle\sup_j d(x,\as_j(t)).
\end{equation}
Since $\as(t)\subset \as_j(t)$, we have $d(x,\as(t))\ge d(x,\as_j(t))$ for every $j$, and hence
$d(x,\as(t))\ge\sup_j d(x,\as_j(t)))$. To prove the opposite inequality, 
we may assume that the right-hand side of \eqref{dxft} is finite. For every $j$ we fix
$y_j\in \as_j(t)$ such that $d(x,y_j)\le d(x,\as_j(t))+\tfrac1j$. As $X$ is compact, there
exists a subsequence $(y_{j_m})$ of $(y_j)$ which converges to a point $y$. Since
$y_{j_m}\in\as_{j_m}(t)$, there exists $k_m\ge j_m$ such that 
$d(f_{k_m}(t),y_{j_m})\le\tfrac1{j_m}$ and $|\vartheta_{k_m}(t)-\vartheta(t)|\le\tfrac1{j_m}$. 
It follows that $f_{k_m}(t)\to y$ and $\vartheta_{k_m}(t)\to \vartheta(t)$ as $m\to\infty$,
hence $y \in \as(t)$ by the very definition of $\as(t)$. Therefore
$$
d(x,\as(t))\le d(x,y)=\lim_{m\to\infty} d(x,y_{j_m})\le \textstyle\sup_j d(x,\as_j(t)),
$$
which concludes the proof of \eqref{dxft}.

On the other hand we have
$
\as_j(t)=\bigcup_{k\ge j}\as^k_j(t)
$,
where
$$
\as^k_j(t):=\{x \in X: d(f_k(t),x)\le\tfrac1j \text{ \,and \,} |\vartheta_k(t)-\vartheta(t)|\le\tfrac1j\},
$$
so that
$$
d(x,\as_j(t))=\inf_{k\ge j} d(x,\as^k_j(t)).
$$
Therefore \eqref{dxft} gives
$$
d(x,\as(t))=\mathop{\smash\sup\vphantom\inf}_j\inf_{k\ge j} d(x,\as^k_j(t)).
$$
Since $t\mapsto f_k(t)$, $t\mapsto \vartheta_k(t)$, and $t\mapsto \vartheta(t)$ 
are measurable, the functions $t\mapsto d(x,\as^k_j(t))$ are measurable for every
$x\in X$, and this concludes the proof of~$(b)$.
\end{dimo}

\subsection{Continuity properties for stress and deformation}
The following proposition shows that, under suitable convexity assumptions, the stress is
left continuous when the crack is left continuous. The same property holds for the deformation and the deformation gradient, 
if some energy terms are strictly convex.

\begin{proposition}
\label{convstresslem}
Let $t \in [0,T]$ be such that $\Gamma(t)=\Gamma^-(t)$ and $\ws$, $-\fs(t)$, and $-\gs(t)$ are convex. Then 
\begin{equation}
\label{convstress}
\partial \ws(\nabla u(s))\weak \partial \ws(\nabla u(t))
\qquad
\text{weakly in }L^{p'}\!(\Om;\msd) \text{ as }s \to t^-.
\end{equation}
If, in addition, $\ws$ and $-\fs(t)$ are strictly convex, then
\begin{eqnarray}
\label{convstressstrong}
&
\nabla u(s)\to\nabla u(t)
\qquad
\text{strongly in }L^p(\Om;\msd) \text{ as }s \to t^-,
\\
\label{convdefstrong}
&
u(s)\to u(t) 
\qquad
\text{strongly in }L^q(\Om;\R^m)\text{ as }s \to t^-,
\\
\label{convbound}
&
u(s)\to u(t) 
\qquad
\text{strongly in }L^{r}(\partial_S \Om,\R^m)\text{ as }s \to t^-.
\end{eqnarray}
\end{proposition}

\begin{dimo}
Since $\Gamma(t)=\Gamma^-(t)$, using
Proposition~\ref{strctgammaprop} it is easy to prove (see  \cite[Section~4.1]{DMFT})  that
for every sequence $s_j \to t^-$ there exists
$\tilde \Gamma$ such that
$$
\Gamma(s_k)\stackrel{\sigma^p}{\to}\tilde \Gamma\qquad\hbox{and}
\qquad \Gamma(t)=\Gamma(0)\cup\tilde \Gamma.
$$
Arguing as in \cite[Sections 4 and 5]{DMFT} we deduce that there exist a subsequence,
not relabelled, and a minimizer $\tilde u$ of $\Eb(t)$ on $AD(\psi(t),\Gamma(t))$ such that
\begin{eqnarray}
&
\label{wsbv}
u(s_j) \weak \tilde u
\quad
\text{weakly in }GSBV^p_q(\Om;\R^m),
\\
&
\label{wstress}
\partial \ws(\nabla u(s_j)) \weak \partial \ws(\nabla \tilde u)
\quad
\text{weakly in }L^{p'}\!(\Om;\msd).
\end{eqnarray}
Since $u(t)$ is a minimizer of $\Eb(t)$ on $AD(\psi(t),\Gamma(t))$, by \eqref{wstress}
the first part of the theorem is proved if we show that the triple
\begin{equation}
\label{triple}
(\partial \ws(\nabla \tilde u), -\partial \fs(t)(\tilde u), -\partial \gs(t)(\tilde u))
\end{equation}
is the same for all minimizers.

\par
To prove this property, we show that for every minimizer $\tilde u$ of $\Eb(t)$ on 
$AD(\psi(t),\Gamma(t))$ the triple \eqref{triple} is a solution 
of a {\it dual problem} which, in our assumptions, admits a unique solution. 
Let us consider the linear space
\begin{equation}
\label{asgammat}
\vub(\Gamma(t)):=\{v \in GSBV^p_q(\Om;\R^m)\,:\,\Sg{0}{v} \tsub \Gamma(t)\},
\end{equation}
where $\Sg{0}{v}$ is defined as in \eqref{Spsi} with $\psi=0$. Note that 
$AD(\psi(t),\Gamma(t))=\psi(t)+\vub(\Gamma(t))=\tilde u+\vub(\Gamma(t))$.
For every 
$$\eta:=(\eta_1,\eta_2,\eta_3) \in 
L^p(\Om;\msd) {\times} 
L^q(\Om;\R^m) {\times}
L^r(\partial_S \Om;\R^m),
$$
we set
$$
\Psi(\eta):=\ws(\eta_1)-\fs(t)(\eta_2)-\gs(t)(\eta_3).
$$
Let us consider
\begin{equation}
\label{phip}
\Phi(\eta):=\inf_{u \in \tilde u +\vub(\Gamma(t))} \Psi\big((\nabla u,u,u)-\eta\big).
\end{equation}
Notice that $\Phi$ is convex, lower semicontinuous, and proper, since 
$$
\Phi(0)=\Psi(\nabla \tilde u,\tilde u,\tilde u)=\Eb(t)(\tilde u)<+\infty.
$$
As a consequence (see \cite{Ro}), we have that
\begin{equation}
\label{dualpb}
\Phi(0)=\Phi^{**}(0)=-\inf \Phi^*,
\end{equation}
where $\Phi^*$ denotes the convex conjugate of $\Phi$, and $\Phi^{**}$ 
denotes the convex conjugate of $\Phi^*$. 
For $\sigma:=(\sigma_1,\sigma_2,\sigma_3) \in 
L^{p'}\!(\Om;\msd) {\times} 
L^{q'}\!(\Om;\R^m) {\times}
L^{r'}\!(\partial_S \Om;\R^m)$
we get
\begin{equation}
\label{phi*}
\begin{array}{c}
\displaystyle
\Phi^*(\sigma)
=\sup_{u,\eta} \left\{
\psc{\sigma}{(\nabla u,u,u)}+\psc{-\sigma}{(\nabla u,u,u)-\eta}
-\Psi((\nabla u,u,u)-\eta)\right\} 
\\ 
\displaystyle
=\psc{\sigma}{(\nabla \tilde u,\tilde u,\tilde u)}+\Psi^*(-\sigma)+
\sup_{v \in \vub(\Gamma(t))}\psc{\sigma}{(\nabla v,v,v)}.
\end{array}
\end{equation}
Since $\vub(\Gamma(t))$ is a linear space, the last term is either zero or $+\infty$, so that we get
\begin{equation}
\label{phi*bis}
\Phi^*(\sigma)=
\begin{cases}
\psc{\sigma}{(\nabla \tilde u,\tilde u,\tilde u)}+\Psi^*(-\sigma)  & \text{if }
\psc{\sigma}{(\nabla v,v,v)}=0 \text{ for all } v \in \vub(\Gamma(t)),
\vspace{2pt}\\
+\infty      & \text{otherwise.}
\end{cases}
\end{equation}
By \eqref{dualpb} we get that
\begin{equation}
\label{basiceq}
-\Psi(\nabla \tilde u,\tilde u,\tilde u)=
\inf_\sigma \left\{ \psc{-\sigma}{(\nabla \tilde u,\tilde u,\tilde u)}+\Psi^*(\sigma):
\psc{\sigma}{(\nabla v,v,v)}=0\text{ for all }v \in \vub(\Gamma(t))\right\}.
\end{equation}
We claim that the triple \eqref{triple} is a solution to the problem on the right-hand side of \eqref{basiceq}. 
To prove this fact we observe that the Euler's equation associated to the minimum problem satisfied by 
$\tilde u$ yields
$$
\psc{\partial \Psi(\nabla \tilde u,\tilde u,\tilde u)}{(\nabla v,v,v)} =0\text{ for all }v \in \vub(\Gamma(t)).
$$
Moreover by duality we have that
\begin{equation*}
\Psi(\nabla \tilde u,\tilde u,\tilde u)+
\Psi^*\big(\partial \Psi(\nabla \tilde u,\tilde u,\tilde u)\big)
=\psc{\partial \Psi(\nabla \tilde u,\tilde u,\tilde u)}{(\nabla \tilde u,\tilde u,\tilde u)},
\end{equation*}
so that we conclude 
\begin{equation*}
-\Psi(\nabla \tilde u,\tilde u,\tilde u)
= -\psc{\partial \Psi(\nabla \tilde u,\tilde u,\tilde u)}{(\nabla \tilde u,\tilde u,\tilde u)}+
\Psi^*\big(\partial \Psi(\nabla \tilde u,\tilde u,\tilde u)\big).
\end{equation*}
Since
$$
\partial \Psi(\nabla \tilde u,\tilde u,\tilde u)=(\partial \ws(\nabla \tilde u), 
-\partial \fs(t)(\tilde u),
-\partial \gs(t)(\tilde u)),
$$ 
this proves our claim.
By the assumption that $\ws$, $-\fs(t)$, and $-\gs(t)$ are convex and $C^1$, we obtain
that
$$
\Psi^*(\sigma)=\ws^*(\sigma_1)+(-\fs(t))^*(\sigma_2)+(-\gs(t))^*(\sigma_3)
$$ 
is strictly convex, so that the problem in the right hand side of \eqref{basiceq} 
admits a unique solution. This proves that \eqref{triple}
is uniquely determined, concluding the proof of the first part of the proposition.

\par
Let us assume now the strict convexity of $\ws$ and $-\fs(t)$. As a consequence $\Eb(t)$ 
is strictly convex, hence $\tilde u=u(t)$. By \eqref{wsbv} this proves
 \eqref{convstressstrong}-\eqref{convbound} with the weak convergence.

Since $\Gamma(t)=\Gamma^-(t)$, we have
that $\Esup(\Gamma(s))\to \Esup(\Gamma(t))$ as $s\to t^-$.
By continuity of the total energy we have that
$$
\lim_{s \to t^-}\Eb(s)(u(s))=\Eb(t)(u(t)),
$$
Using the lower semicontinuity of all terms in \eqref{elener}  we get
\begin{equation}
\label{convws}
\lim_{s \to t^-}\ws(\nabla u(s))=\ws(\nabla u(t))
\quad\hbox{and}\quad
\lim_{s \to t^-}\fs(s)(u(s))=\fs(t)(u(t)).
\end{equation}
Then the strong convergence in \eqref{convstressstrong} and \eqref{convdefstrong} can be derived from
\eqref{convws}, using a general argument which allows to deduce strong convergence 
from weak convergence and convergence of strictly convex energies (see \cite{Bre2}). 
The strong convergence in \eqref{convbound} follows from (3.27) in~\cite{DMFT}. 
\end{dimo}

{
}
\end{document}